%% file: ruprecht_contrib.tex
\documentclass[runningheads]{lncse}

\usepackage{amsmath}
\usepackage{amsfonts}
\usepackage{epsfig}
\usepackage{graphics} 

\begin{document}

\title{Convergence of Parareal for the Navier-Stokes equations depending on the Reynolds number}

\titlerunning{Parareal for the Navier-Stokes equations}

\author{ Johannes Steiner\inst{1} \and Daniel Ruprecht\inst{1} \and Robert Speck\inst{2,1} \and Rolf Krause\inst{1}}

\authorrunning{Steiner, Ruprecht, Speck, Krause}   

\institute{
Institute of Computational Science, Universit{\`a} della Svizzera italiana, Via Giuseppe Buffi 13, CH-6900 Lugano, Switzerland, {\tt \{johannes.steiner, daniel.ruprecht, rolf.krause\}@usi.ch}
\and J\"ulich Supercomputing Centre, Forschungszentrum J\"ulich, J\"ulich, Germany, {\tt r.speck@fz-juelich.de}
}

\maketitle

\begin{abstract}
The paper presents first a linear stability analysis for the time-parallel Parareal method, using an IMEX Euler as coarse and a Runge-Kutta-3 method as fine propagator,
confirming that dominant imaginary eigenvalues negatively affect Parareal's convergence.
This suggests that when Parareal is applied to the nonlinear Navier-Stokes equations, problems for small viscosities could arise.
Numerical results for a driven cavity benchmark are presented, confirming that Parareal's convergence can indeed deteriorate as viscosity decreases and the flow becomes increasingly dominated by convection. The effect is found to strongly depend on the spatial resolution.
\end{abstract}

\section{Introduction}
As core counts in modern supercomputers continue to grow, parallel algorithms are required that can provide concurrency beyond existing approaches parallelizing in space. In particular, algorithms that parallelize in time "along the steps" have attracted noticeable interest. 
Probably the most widely studied algorithm of this type is Parareal~\cite{ruprecht_contrib_LionsEtAl2001}, but other important methods exist as well, for example PITA~\cite{ruprecht_contrib_FarhatEtAl2003} or PFASST~\cite{ruprecht_contrib_EmmettMinion2012}.\par
The applicability of Parareal to the Navier-Stokes equations has been studied in~\cite{ruprecht_contrib_FischerEtAl2005}, where it is shown that Parareal can solve the initial value problem arising from a Finite Element discretization of the Navier-Stokes equations for a Reynolds number of 200 as well as from a Spectral Element discretization for a problem with Reynolds number 7,500.
A non-Newtonian problem is studied in~\cite{ruprecht_contrib_Celledoni2009}.  
In~\cite{ruprecht_contrib_Trindade2004,ruprecht_contrib_Trindade2006}, Parareal is combined with parallelization in space and setups with Reynolds numbers up to 1,000 are investigated. 
While it is confirmed that Parareal can successfully be applied to flow simulations, the attempt to demonstrate its potential to provide speedup beyond the saturation of the spatial parallelization was inconclusive, as either the pure time or pure space parallel approach provided minimum runtimes. 
A successful demonstration that Parareal can speed up simulations after the spatial parallelization has saturated can be found in~\cite{ruprecht_contrib_CroceEtAl2012}, where Parareal is used to simulate a driven cavity flow in a cube with a Reynolds number of 1,000.
The performance of PFASST for a particle-based discretization of the Navier-Stokes equations on $\mathcal{O}(100,000)$ cores is studied in~\cite{ruprecht_contrib_SpeckEtAl2012}.\par
It has been noted in multiple works that Parareal as well as PITA have stability issues for convection-dominated problems, see~\cite{ruprecht_contrib_Bal2005,ruprecht_contrib_FarhatEtAl2003,ruprecht_contrib_GanderVandewalle2007_SISC,ruprecht_contrib_RuprechtKrause2012,ruprecht_contrib_StaffRonquist2005}.
This suggests that Parareal will at some point cease to converge properly for the Navier-Stokes equations if the Reynolds number increases and the problem becomes more and more dominated by advection.
This paper discusses results from linear stability analysis and presents a numerical study for two-dimensional driven cavity flow of how the convergence of Parareal is affected as viscosity decreases.

\section{Parareal}
Parareal is a method to introduce concurrency in the solution of initial value problems
\begin{equation}
  \label{ruprecht_contrib_eq:ivp}
    u_{t} = f(u(t),t), \quad u(0) = u_{0}, \quad 0 \leq t \leq T.
\end{equation}
It relies on the introduction of two classical one-step time integration methods, one computationally expensive and of high accuracy (denoted by $\mathcal{F}$) and one computationally cheap method of lower accuracy (denoted by $\mathcal{G}$). 
The former is commonly referred to as the "fine propagator", the latter as the "coarse propagator". 
Denote by $U_n$ the numerical approximation of the exact solution $u$ of~\eqref{ruprecht_contrib_eq:ivp} at some point in time $t_n$. 
Further, denote as
\begin{equation}
    \label{ruprecht_contrib_eq:serial_step}
    U_{n+1} = \mathcal{F}_{\delta t}(U_n)
\end{equation}
the result obtained by integrating from an initial value $U_n$ given at a time $t_{n}$ forward in time to a time $t_{n+1}$ using a time-step $\delta t$ and the method indicated by $\mathcal{F}$.
For a decomposition of $[0,T]$ into $N$ so-called time-slices $[t_n, t_{n+1}]$, $n=0,\ldots,N-1$, solving~\eqref{ruprecht_contrib_eq:serial_step} time-slice after time-slice corresponds to classical time-marching, running the fine method in serial from $t_0 = 0$ to $t_{N} = T$.
Instead, Parareal approximately computes the values $U_{n}$ by means of the iteration
\begin{equation}
	\label{ruprecht_contrib_eq:parareal}
	U^{k+1}_{n+1} = \mathcal{G}_{\Delta t}(U^{k+1}_{n}) + \mathcal{F}_{\delta t}(U^{k}_{n}) - \mathcal{G}_{\Delta t}(U^{k}_{n})
\end{equation}
were $k$ denotes the iteration counter.
For $k \to N$, iteration~\eqref{ruprecht_contrib_eq:parareal} converges towards the serial fine solution, that is $U^{k}_{n} \to U_{n}$.
Once values $U^{k}_{n}$ are known, the evaluation of the computationally expensive terms $\mathcal{F}(U_n^{k})$ in~\eqref{ruprecht_contrib_eq:parareal} can be done in parallel on $N$ processors.
Then, a correction is propagated serially by evaluating the terms $\mathcal{G}_{\Delta t}(U^{k+1}_{n})$ and computing $U^{k+1}_{n+1}$.
We refer to e.g.~\cite{ruprecht_contrib_RuprechtKrause2012} for a more in-depth presentation of the algorithm.
The speedup achievable by Parareal concurrently computing the solution on $N$ time-intervals assigned to $N$ processors is bounded by
\begin{equation}
	\label{ruprecht_contrib_eq:speedup_bounds}
  s(N) \leq \min\left\{ \frac{N}{N_{\rm it}}, \frac{ C_{\mathcal{F}} }{ C_{\mathcal{G}} } \right\}
\end{equation}
where $N_{\rm it}$ is the number of iterations performed and $C_{\mathcal{F}}$, $C_{\mathcal{G}}$ denote the time required to evaluate $\mathcal{F}_{\delta t}$ and $\mathcal{G}_{\Delta t}$ respectively, see again e.g.~\cite{ruprecht_contrib_RuprechtKrause2012}.
Note that the two bounds are competing in the sense that using a coarser and cheaper method for $\mathcal{G}$ will usually improve the second bound but might cause Parareal to require more iterations to converge, thereby reducing the first bound.
In contrast, a more accurate and more expensive $\mathcal{G}$ will likely reduce the iteration number but also reduce the coarse-to-fine runtime ratio $\frac{C_{\mathcal{F}}}{C_{\mathcal{G}}}$.

\section{Linear stability analysis}
In order to illustrate Parareal's stability properties, we apply it to the test equation
\begin{equation}
	\label{eq:test}
	y'(t) = \lambda_{\rm Re} y(t) + i \lambda_{\rm Im} y(t), \quad y(0) = 0, \quad 0 \leq t \leq T.
\end{equation}
A linear stability analysis of this kind was first done in~\cite{ruprecht_contrib_StaffRonquist2005}, using RadauIIA methods for both $\mathcal{F}$ and $\mathcal{G}$.
Here, in line with the numerical examples presented in Section~\ref{sec:results}, the stability analysis is done for an implicit-explicit Euler method for $\mathcal{G}$ and an explicit Runge-Kutta-3 method for $\mathcal{F}$ with five time steps of $\mathcal{F}$ per two time steps of $\mathcal{G}$.
The IMEX scheme treats the real part ("diffusion") implicitly and the imaginary term ("convection") explicitly.
Further, $N = 15$ concurrent time slices are used and a time step $\Delta t = 1.0$ for $\mathcal{G}$, so that $T = 15$.\par
\begin{figure}[t]
	\centering
	\begin{minipage}{0.325\textwidth}
		\centering
		(a) IMEX Euler
		\includegraphics[width=\textwidth]{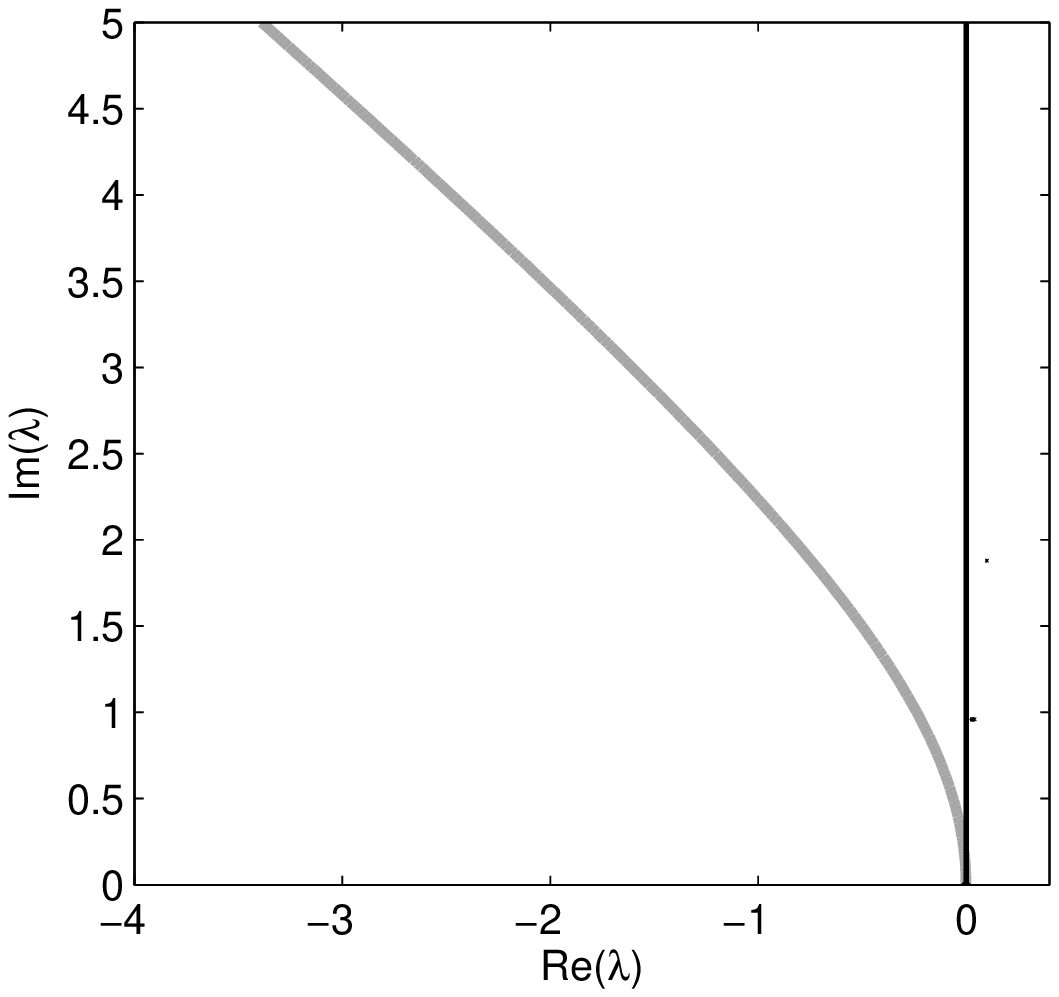}\vspace{1em}
	\end{minipage}
	\begin{minipage}{0.325\textwidth}
		\centering
		(b) RK-3
		\includegraphics[width=\textwidth]{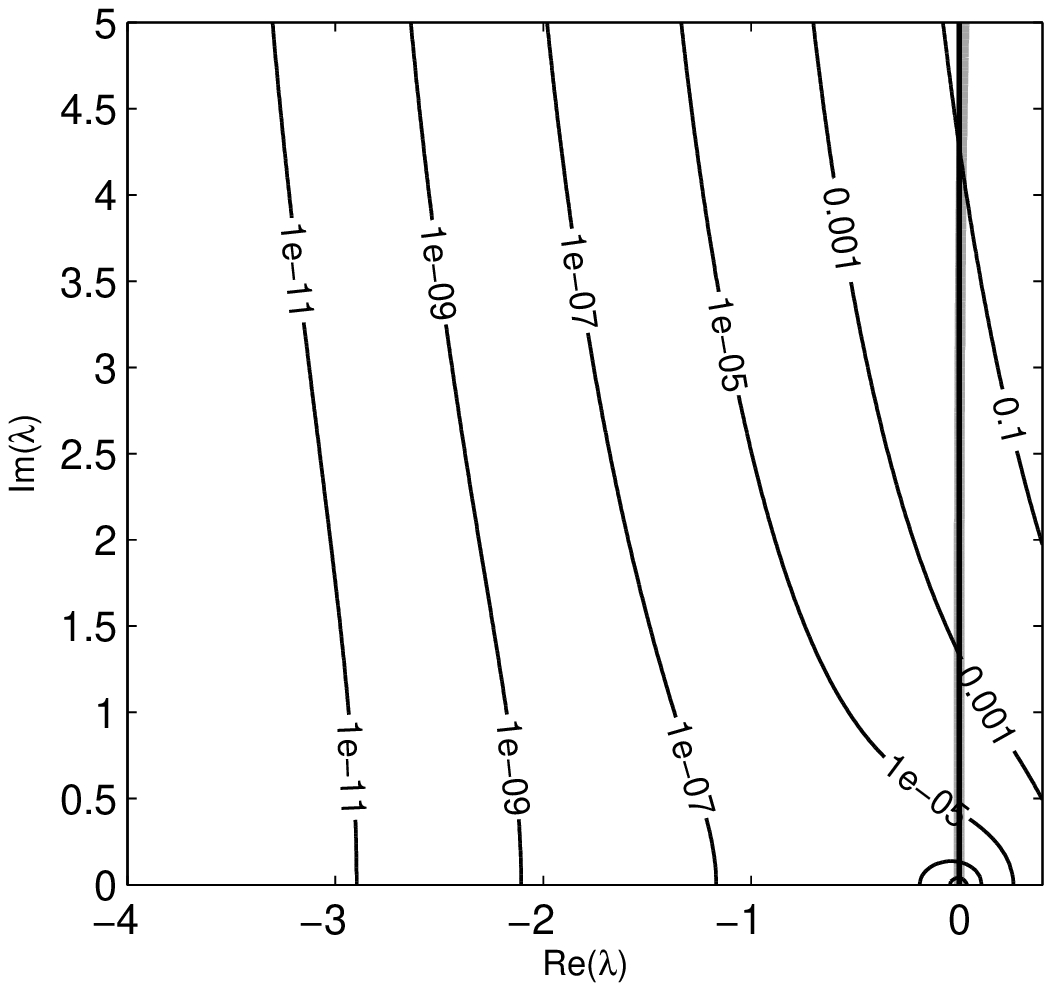}\vspace{1em}
	\end{minipage}
	\begin{minipage}{0.325\textwidth}
		\centering
		(c) Parareal(1)
		\includegraphics[width=\textwidth]{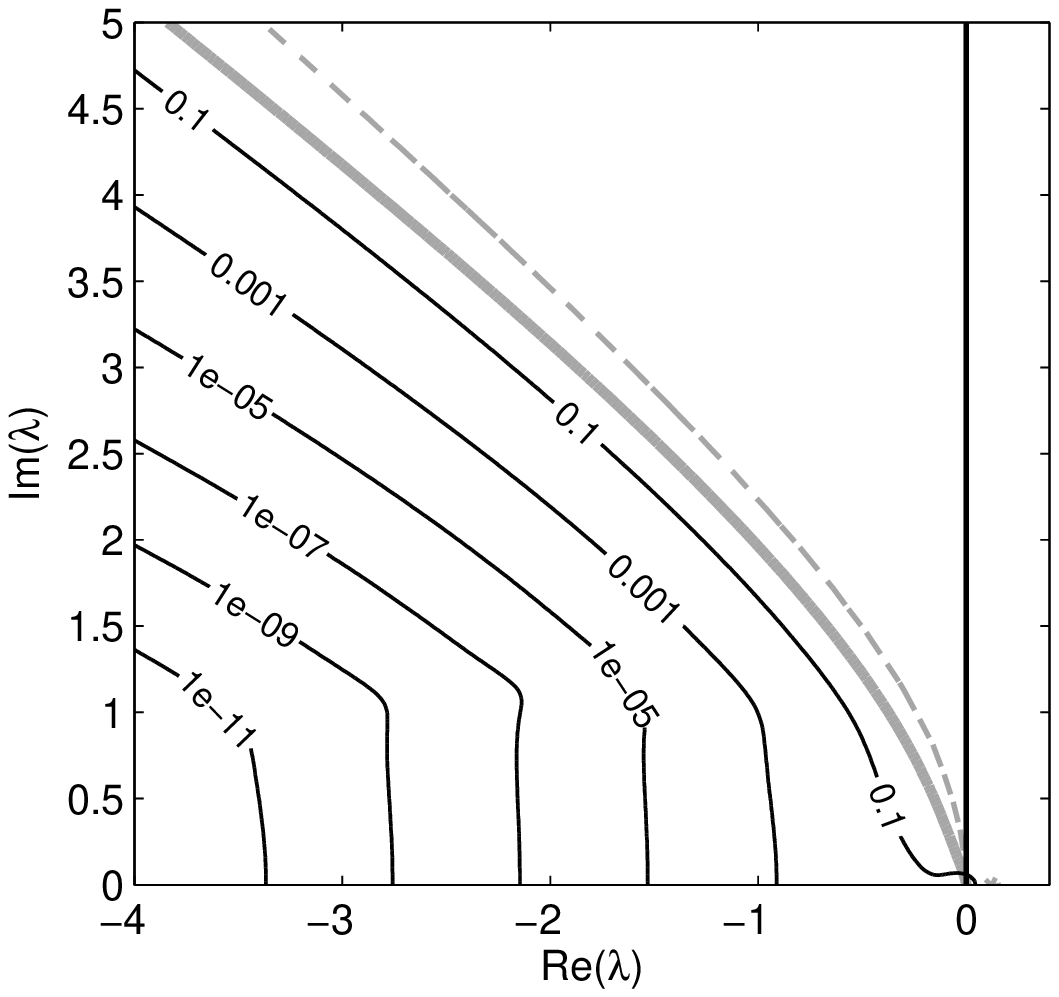}\vspace{1em}
	\end{minipage}
	\begin{minipage}{0.325\textwidth}
		\centering
		(d) Parareal(4)
		\includegraphics[width=\textwidth]{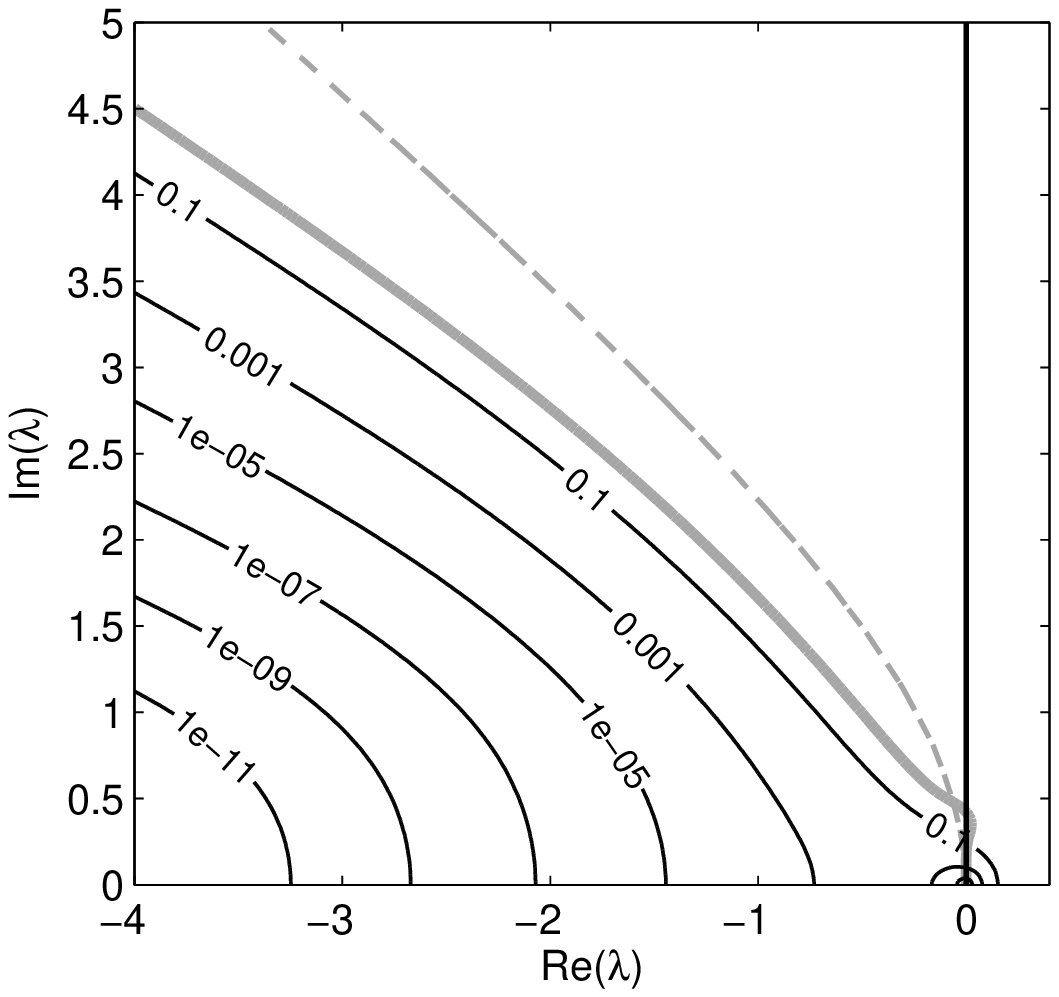}
	\end{minipage}
	\begin{minipage}{0.325\textwidth}
		\centering
		(e) Parareal(8)
		\includegraphics[width=\textwidth]{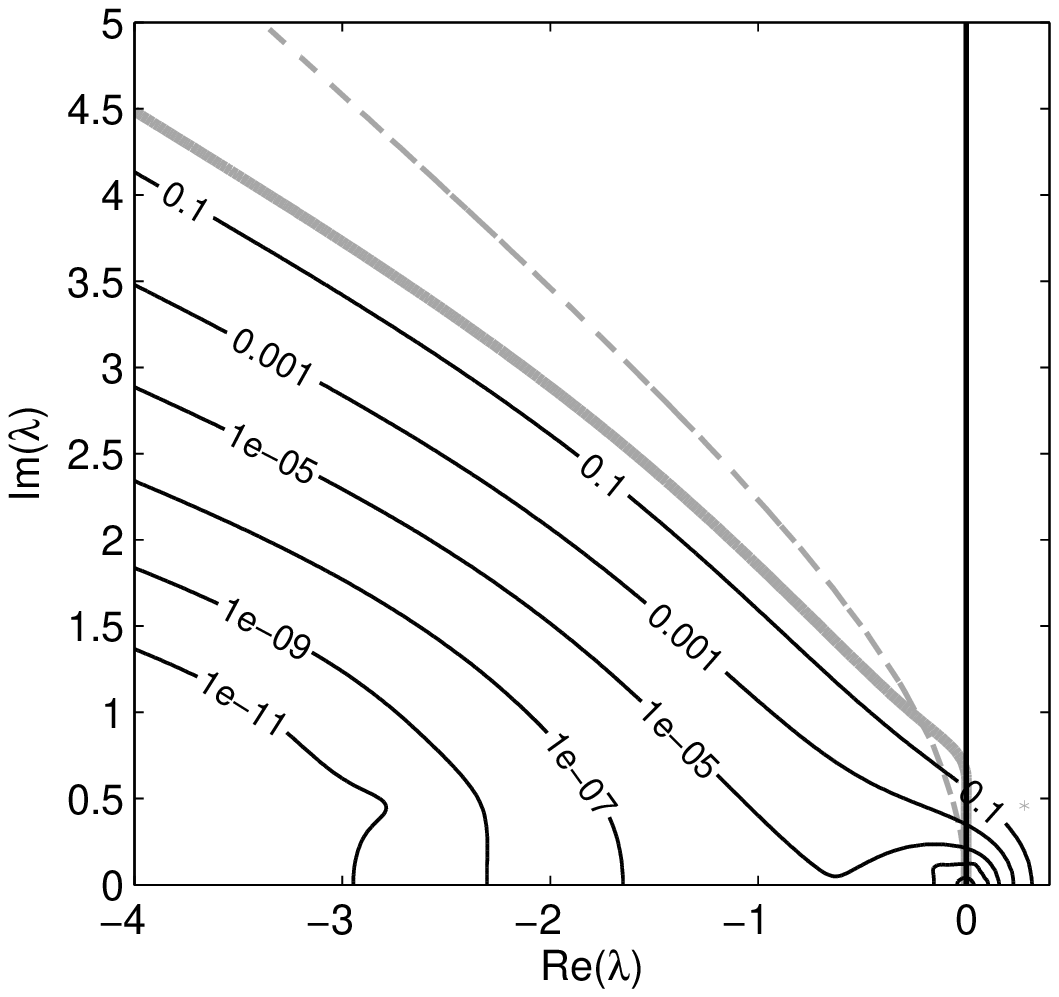}
	\end{minipage}
	\begin{minipage}{0.325\textwidth}
		\centering
		(f) Parareal(12)
		\includegraphics[width=\textwidth]{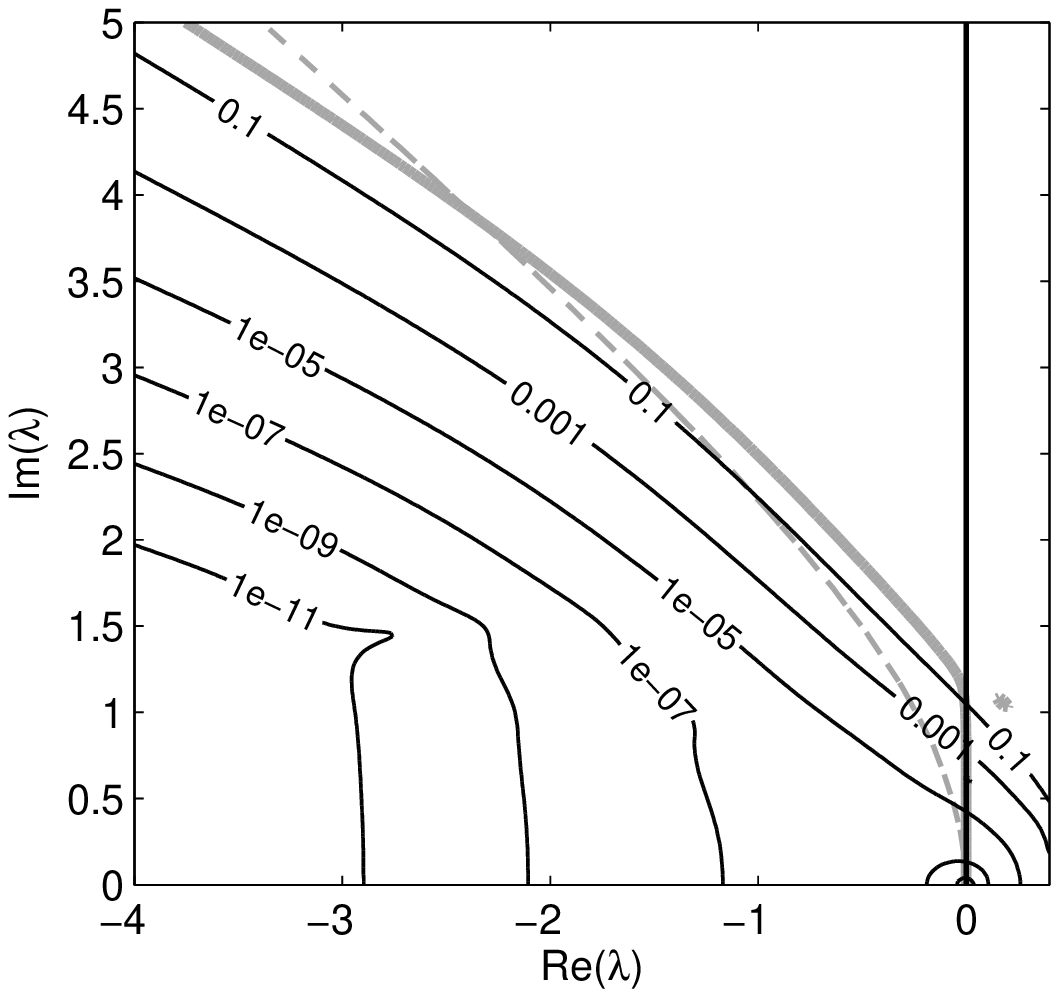}
	\end{minipage}
	\caption{Stability and accuracy of Parareal using an implicit-explicit Euler for $\mathcal{G}$, a RK3 method for $\mathcal{F}$, $N = 15$ time slices and a ratio of $s=5/2$ fine to coarse steps in each time-slice. The thick gray line indicates where the amplification factor becomes greater than one. The black lines indicate error levels. Note that in (a) no black lines are visible because the error never drops below $10^{-1}$. Note also that $s=5/2$ means the fine scheme in serial performs five steps per time-slice and the coarse scheme two, so that (a) and (b) are not identical to the stability function of the respective method with only a single time-step. Figures (c) -- (f) show the stability domain for Parareal with $N_{\rm it} = 1, 4, 8, 12$ iterations. For comparison, the stability region of $\mathcal{G}$ is also sketched again as a thin dashed gray line.}\label{ruprecht_contrib_fig:stability}
\end{figure}

Figure~\ref{ruprecht_contrib_fig:stability} shows the resulting stability domains and isolines of accuracy for the coarse method run serially (a), the fine method run serially (b), and for Parareal with different numbers of iteration (c)--(f).
For $N_{\rm it} = N = 15$, the solution from Parareal is identical to the one provided by $\mathcal{F}$ and thus the stability domains also coincide (not shown).
As can be expected because of the stability constraint arising from the explicitly treated imaginary term, the IMEX method used for $\mathcal{G}$ becomes unstable if the imaginary part of $\lambda$ becomes too dominant.
Parareal however ceases to be stable even before reaching the stability limit of the coarse propagator.
The analysis confirms again that for problems with imaginary eigenvalues, Parareal can develop instabilities although both $\mathcal{F}$ and $\mathcal{G}$ are stable.
Furthermore, the stability domain of Parareal shrinks from $N_{\rm it} = 1$ to $N_{\rm it} = 4$ and $N_{\rm it} = 8$ before expanding again for $N_{\rm it} = 12$.
Note also that for a fixed number of iterations, Parareal becomes less accurate as $\lambda_{\rm Im}$ increases (in contrast to the serial fine method), corresponding to reduced rates of convergence.
This means that achieving the accuracy of the underlying fine method will require more iterations for problems with larger imaginary eigenvalues, therefore reducing the speedup achievable by Parareal, cf. the estimate~\eqref{ruprecht_contrib_eq:speedup_bounds}.
Eventually, as convergence becomes too slow, Parareal will no longer be able to achieve speedup at all and will no longer be useful.
The mathematical explanation for this behavior is a growing term in the error estimate for Parareal for imaginary eigenvalues that is only compensated for as the iteration number approaches the number of time-slices, see the analysis in~\cite{ruprecht_contrib_GanderVandewalle2007_SISC}.

\section{Numerical results for driven cavity flow}\label{sec:results}
In order to investigate if and how the results from the linear stability analysis carry over to the fully nonlinear case, we solve now the non-dimensional, nonlinear, incompressible Navier-Stokes equations in two dimensions
\begin{align}
	\label{ruprecht_contrib_eq:n_s_e}
        \mathbf{u}_{t} + \mathbf{u} \cdot \nabla \mathbf{u} + \nabla p &= \nu \Delta \mathbf{u} \\
        \nabla \cdot \mathbf{u} &= 0
\end{align}
on a square $[0,1]^{2}$.
A method-of-lines approach is used to first discretize in space. For the spatial discretization a finite volume method based on a vertex centered scheme is used. 
On an unstructured or not necessarily structured triangle mesh, control volumes are constructed via a dual mesh.
This leads to a non-staggered scheme of velocity and pressure. Therefore, a stabilization based on upwind differences and an incremental version of the Chorin-Temam method for the pressure is used~\cite{ruprecht_contrib_versteeg2007introduction}.
\begin{figure}[t]
	\centering
	\begin{minipage}{0.495\textwidth}
		\centering
		(a) $N_x = 8$
		\includegraphics[width=\textwidth]{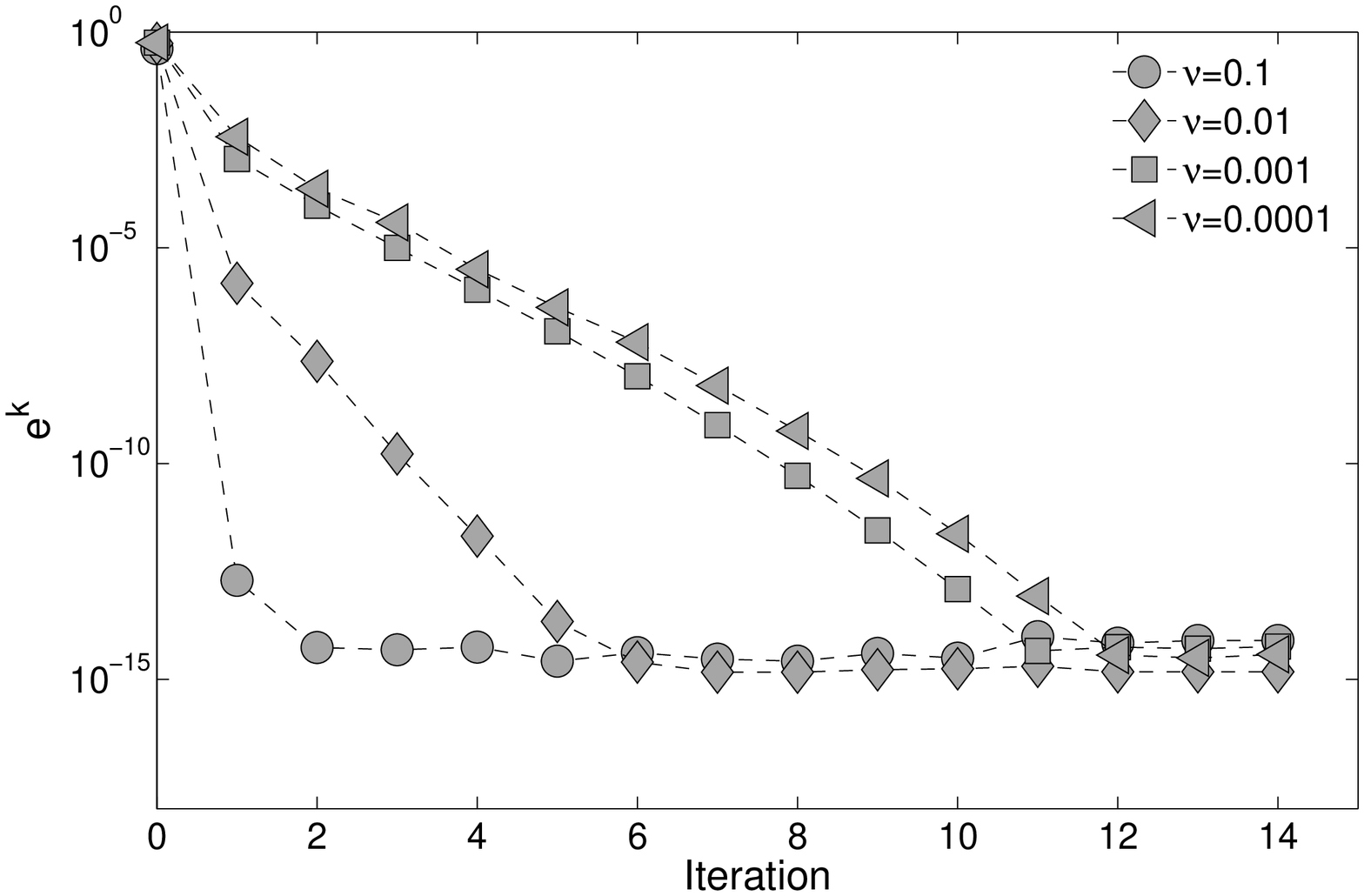}\vspace{0.5em}
	\end{minipage}
	\begin{minipage}{0.495\textwidth}
		\centering
		(b) $N_x = 16$
		\includegraphics[width=\textwidth]{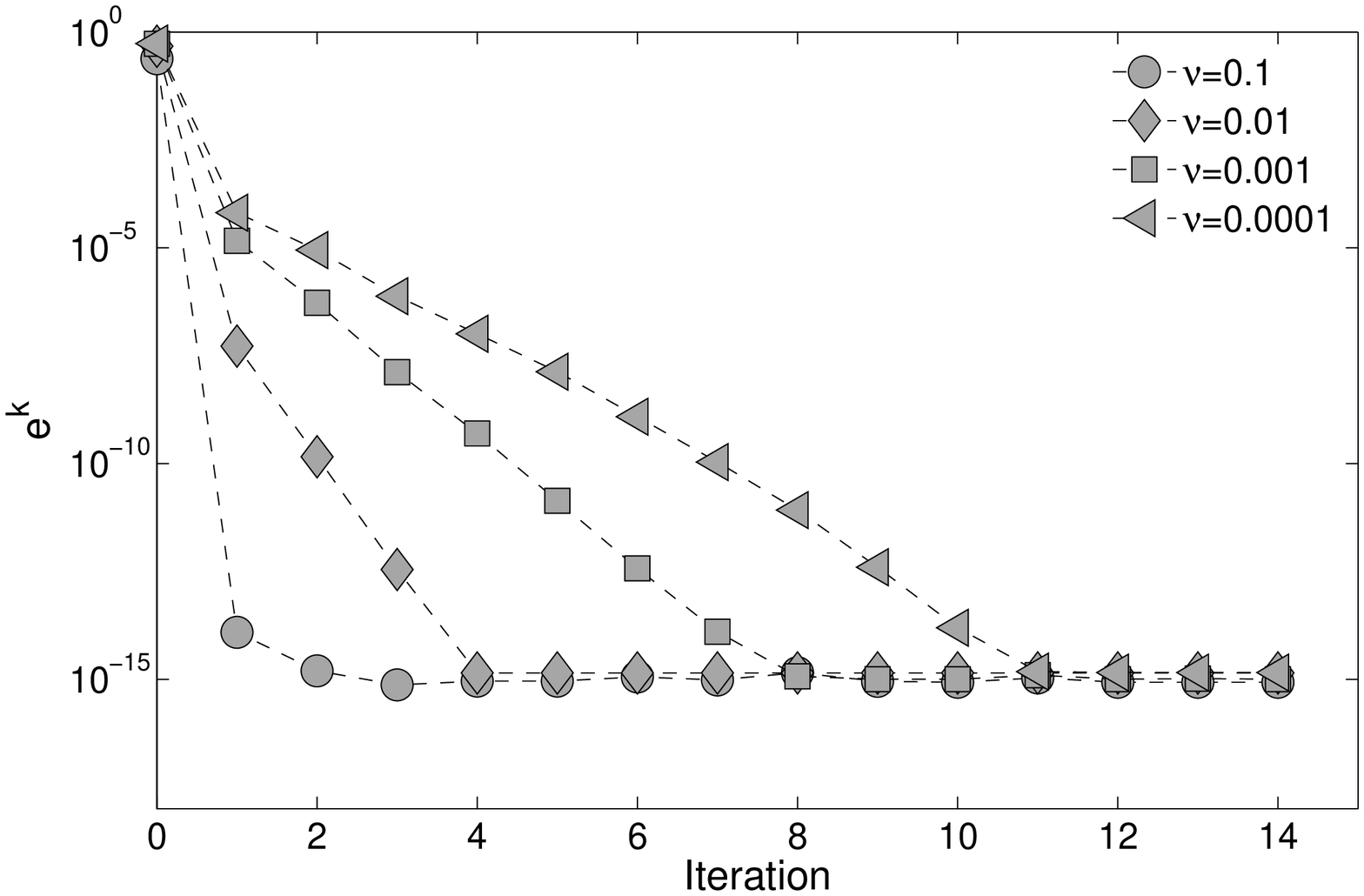}\vspace{0.5em}
	\end{minipage}
	\begin{minipage}{0.495\textwidth}
		\centering
		(c) $N_x = 32$
		\includegraphics[width=\textwidth]{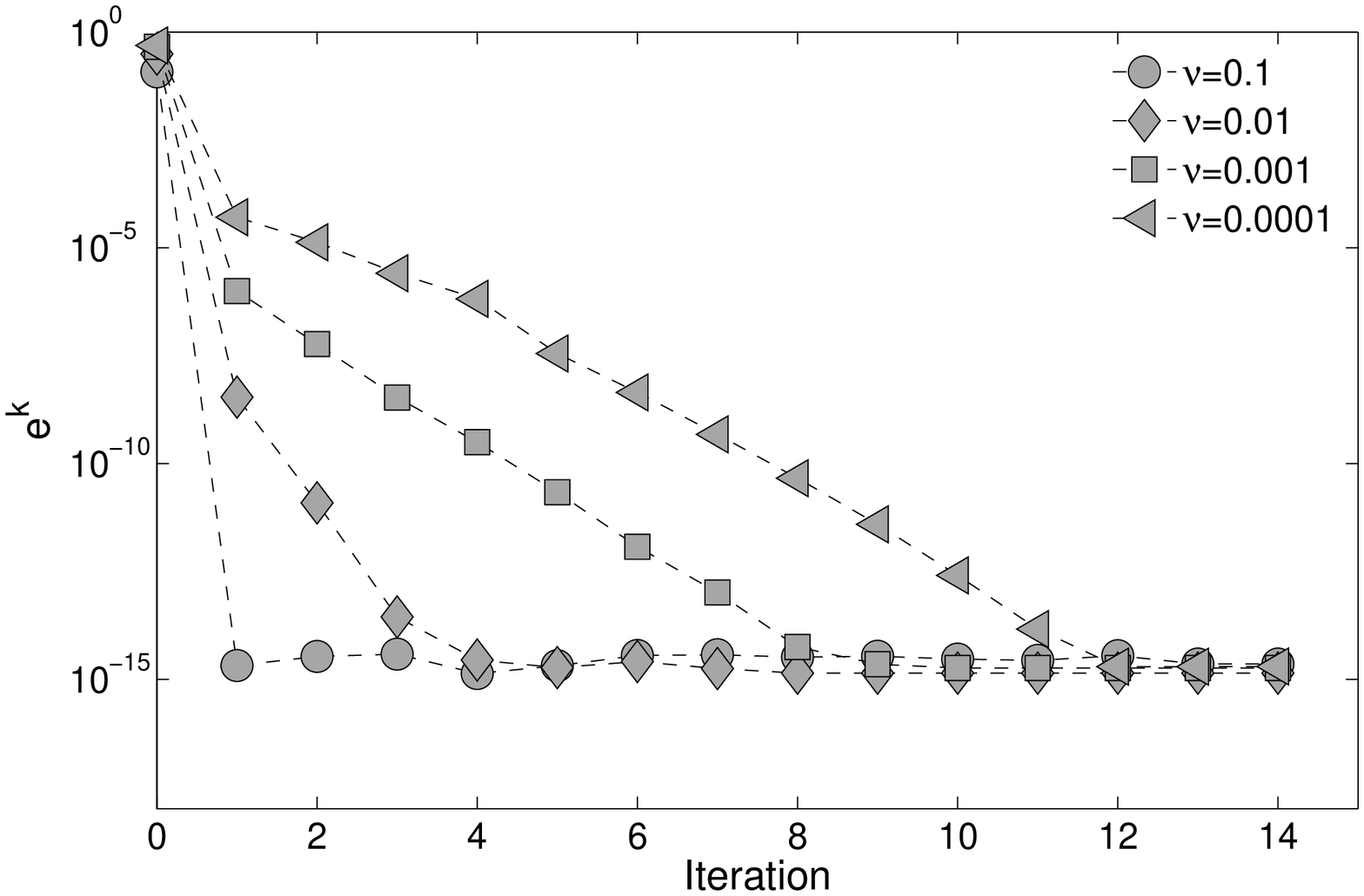}
	\end{minipage}
	\begin{minipage}{0.495\textwidth}
		\centering
		(d) $N_x = 64$
		\includegraphics[width=\textwidth]{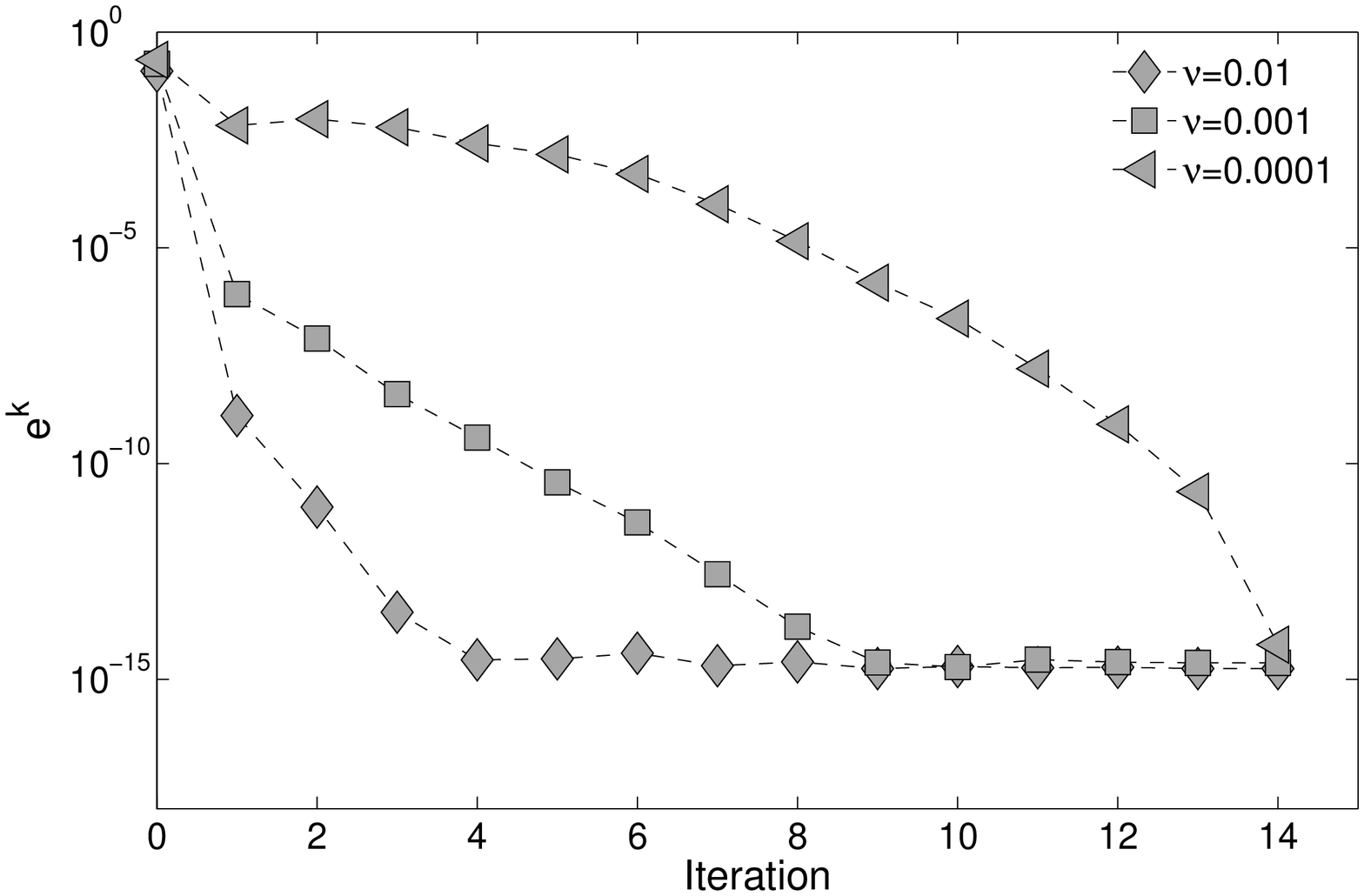}
	\end{minipage}
	\caption{Convergence of Parareal against the serial fine solution for $\Delta t=1/200$, different numbers of mesh points $N_x$ and different values for the viscosity $\nu$.}\label{ruprecht_contrib_fig:para_conv}
\end{figure}
Parareal is then employed to solve the resulting initial value problem until a final time $T=15$ with $N=15$ time-slices.
As in the stability analysis above, $\mathcal{G}$ is an implicit-explicit Euler method while $\mathcal{F}$ is an explicit Runge-Kutta-3 method.
The time-step for the coarse method is $\Delta t = 1/200$, for the fine method $\delta t = 1/500$, reproducing a rate of $s=5/2$ fine per coarse steps.
Although the driven cavity setup is probably not the most ideal here, since, depending on the viscosity, the solution settles into a steady state rather quickly, its wide use and comparative simplicity still make for a good first test case.
Further tests for a more complex vortex shedding setups are currently ongoing.\par
Figure~\ref{ruprecht_contrib_fig:para_conv} shows the convergence of Parareal against the solution provided by running $\mathcal{F}$ in serial. 
Shown is the maximum of the relative error at the end of all time-slices, that is
\begin{equation}
	e^{k}  := \max_{n=1,\ldots,N} \left\| U^{k}_{n} - U_{n} \right\|_{\infty}/\left\| U_{n} \right\|_{\infty}
\end{equation}	
\begin{figure}[t]
	\centering
	\begin{minipage}{0.495\textwidth}
		\centering
		(a) $N_x = 8$
		\includegraphics[width=\textwidth]{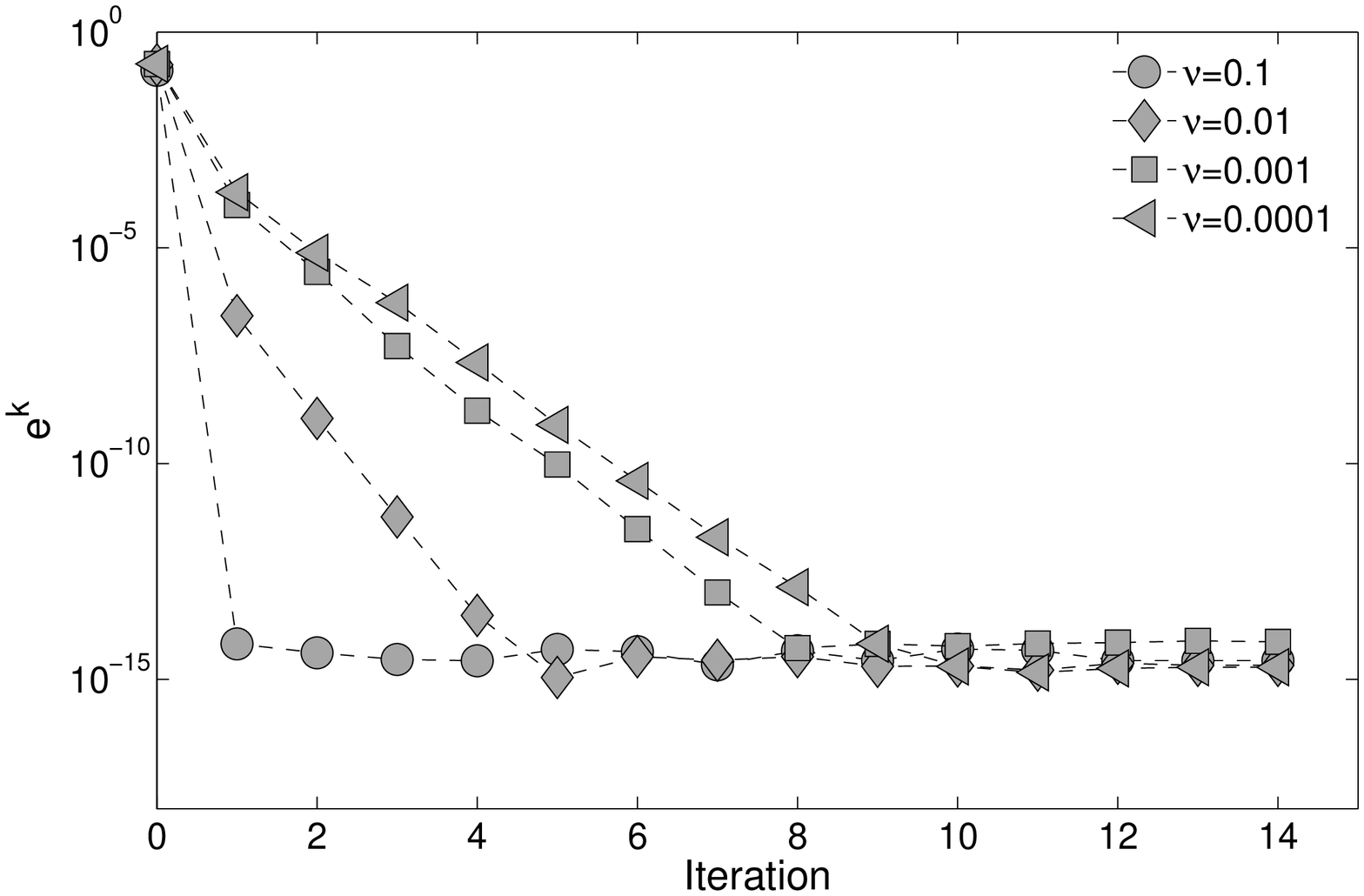}\vspace{0.5em}
	\end{minipage}
	\begin{minipage}{0.495\textwidth}
		\centering
		(b) $N_x = 16$
		\includegraphics[width=\textwidth]{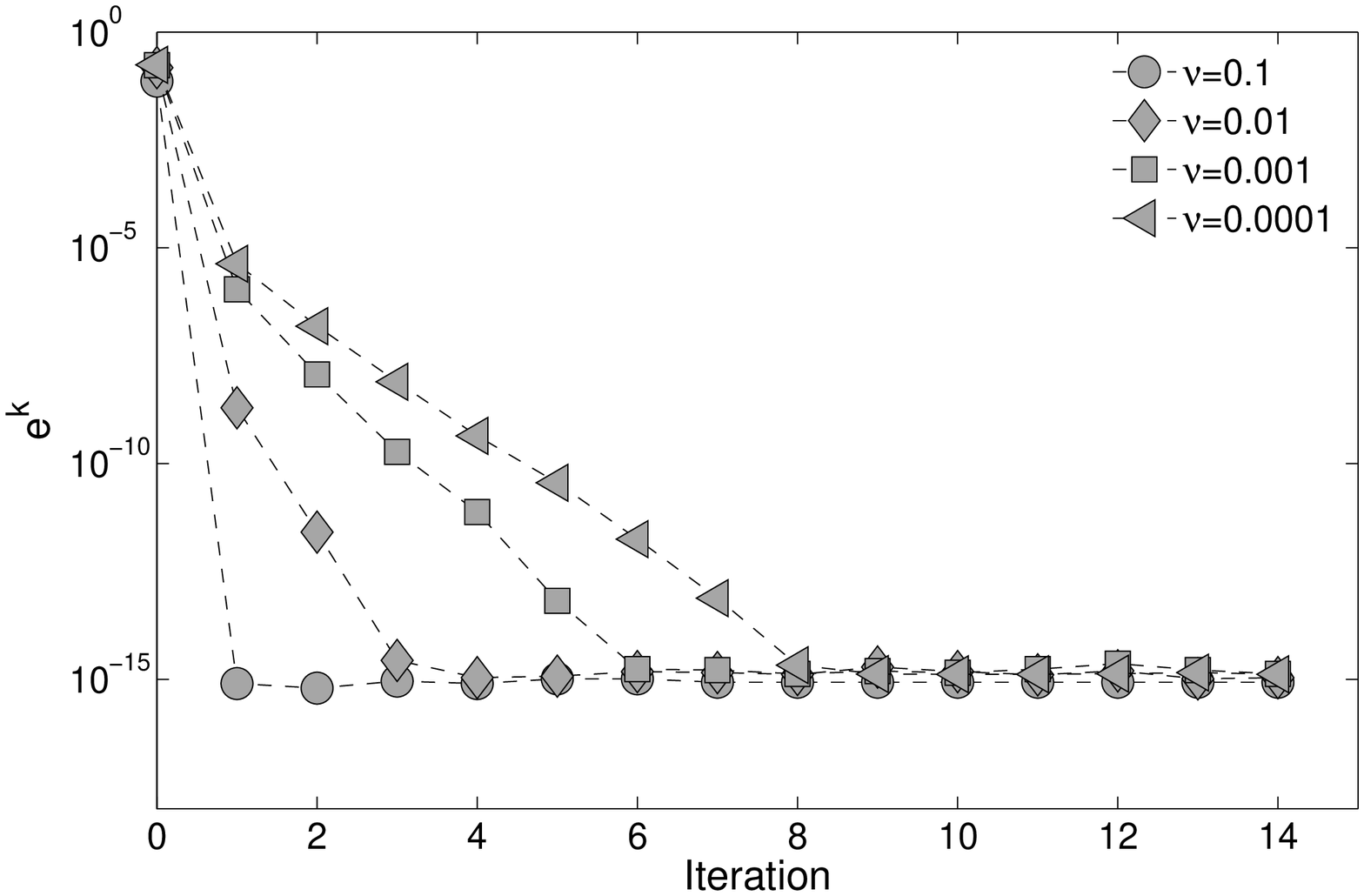}\vspace{0.5em}
	\end{minipage}
	\begin{minipage}{0.495\textwidth}
		\centering
		(c) $N_x = 32$
		\includegraphics[width=\textwidth]{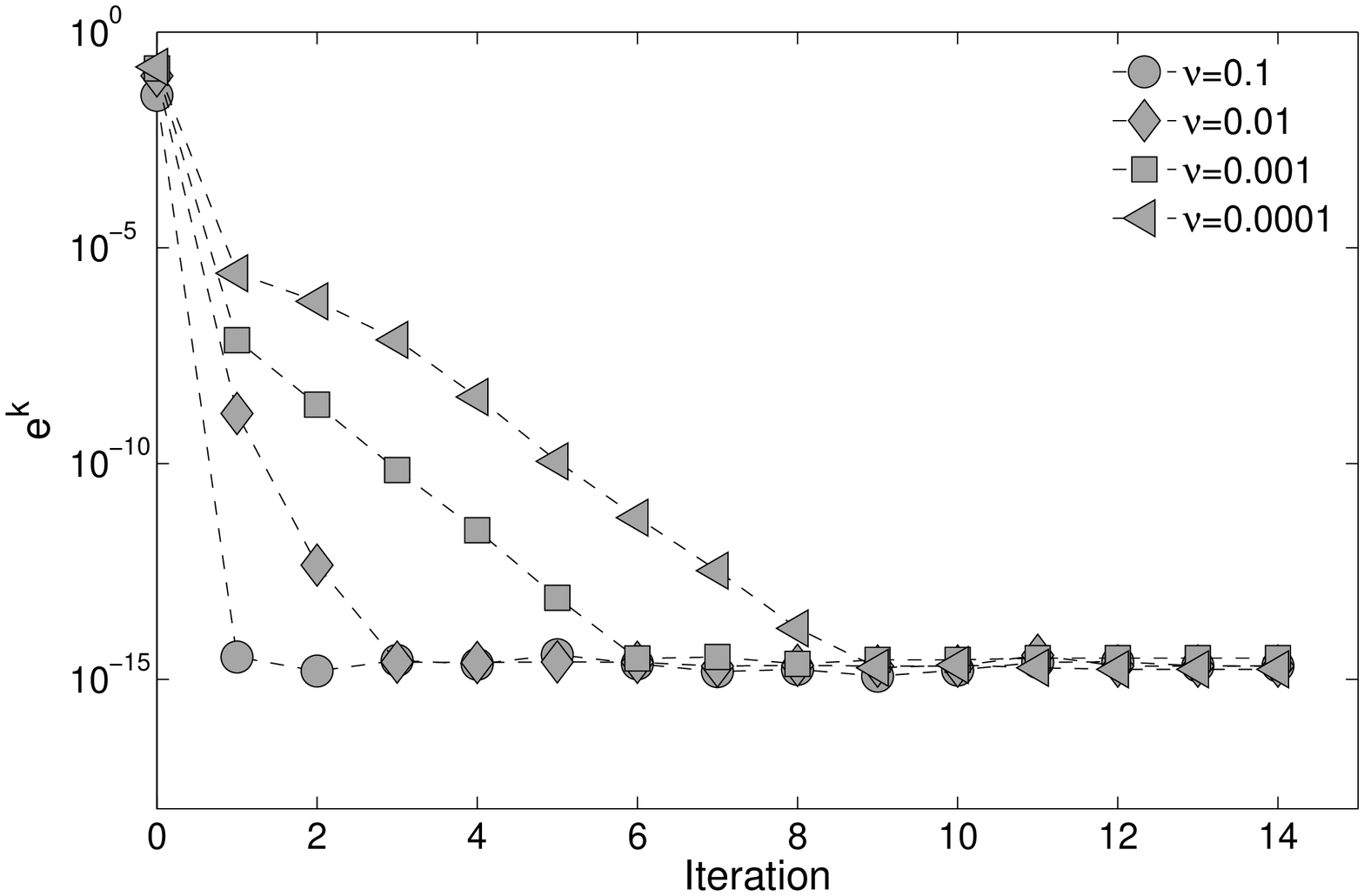}
	\end{minipage}
	\begin{minipage}{0.495\textwidth}
		\centering
		(d) $N_x = 64$
		\includegraphics[width=\textwidth]{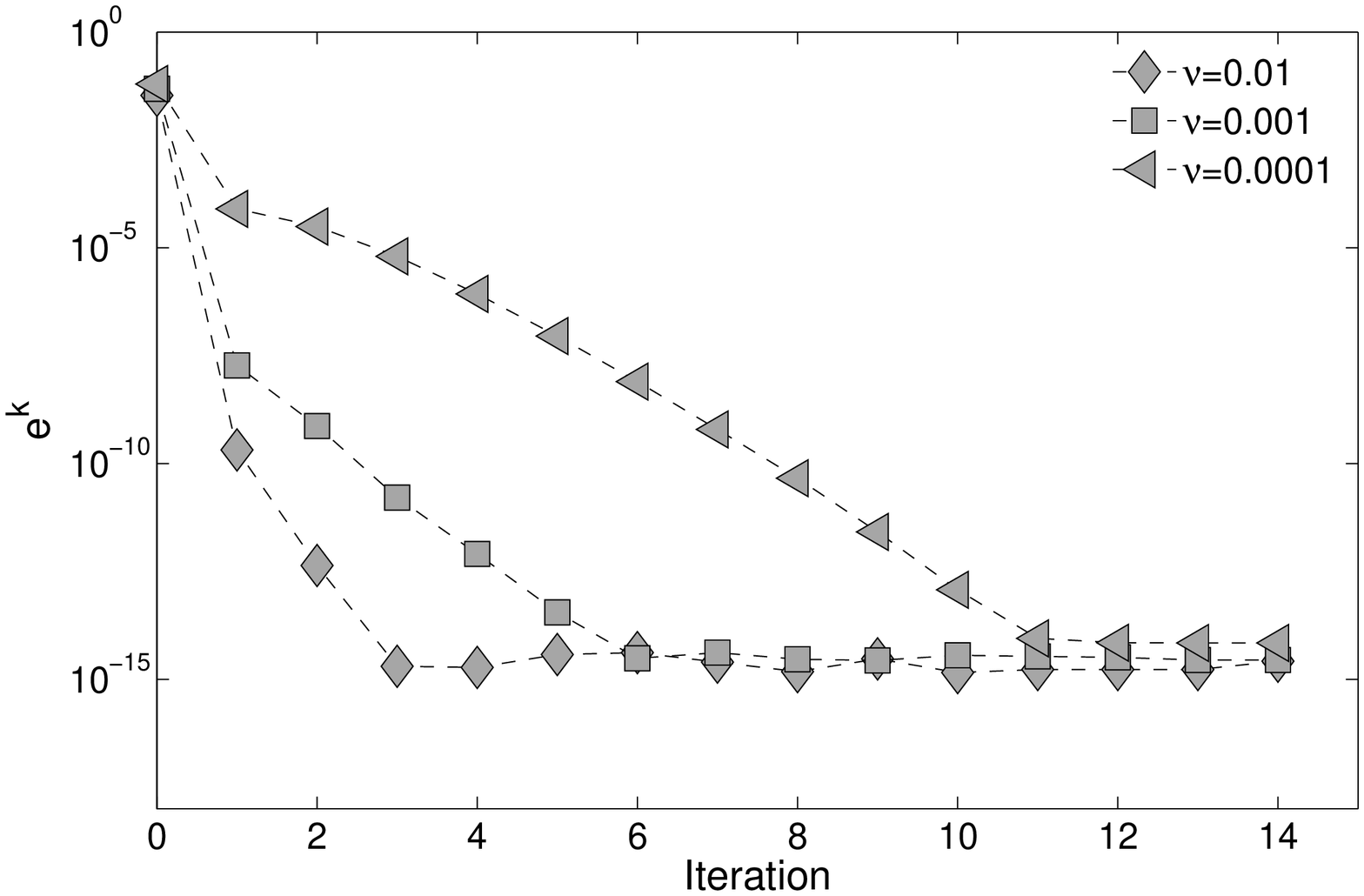}
	\end{minipage}
	\caption{Convergence of Parareal against the serial fine solution for $\Delta t = 1/400$, different numbers of mesh points $N_x$ and different values for the viscosity $\nu$.}\label{ruprecht_contrib_fig:para_conv2}
\end{figure}
where $U^{k}_{n}$ is the solution at $t_{n}$ provided by Parareal after $k$ iterations and $U_{n}$ the solution provided by running $\mathcal{F}$ in serial.
The spatial discretization uses values of $N_x=8$, $16$, $32$, $64$ and the viscosity parameter is set to $\nu=10^{-1}$, $10^{-2}$, $10^{-3}$, $10^{-4}$.
For $N_x=64$ and $\nu = 10^{-1}$ no values are shown, because here the explicit RK3 method used for $\mathcal{F}$ started to show stability problems.
On all meshes, the convergence of Parareal deteriorates as $\nu$ becomes smaller and this effect is much more pronounced for finer spatial resolutions, where the mesh is able to better resolve the features of the more convection dominated flow.
On the finest mesh, there is a clear transition between $\nu=10^{-3}$, for which Parareal still converges reasonably well, and $\nu = 10^{-4}$, where the method first stalls for several iterations before slowly starting to converge.
Requiring a number of iterations close to the number of time-slices means that only marginal speedup is possible from Parareal, because the first bound in~\eqref{ruprecht_contrib_eq:speedup_bounds} becomes very small.
Note also that the still reasonable convergence of Parareal for very low viscosity on a very coarse spatial mesh is not of great practical interest, as the provided solution will be strongly under-resolved.
Figure~\ref{ruprecht_contrib_fig:para_conv2} shows again the convergence of Parareal for a decreased coarse time-step size $\Delta t = 1/400$. 
As can be seen, reducing the coarse time-step again improves convergence and allows Parareal to converge in fewer iterations.
However, it reduces the second speedup bound in~\eqref{ruprecht_contrib_eq:speedup_bounds} and thus will also at some point prevent Parareal from achieving speedup.
Therefore, the reduced convergence speed of Parareal for small viscosities either necessitates a small time-step in the coarse method or a large number of iterations and both choices significantly reduce the achievable speedup.
A possible remedy could be the application of stabilization techniques as discussed in \cite{ruprecht_contrib_CortialFarhat2009,ruprecht_contrib_FarhatCortial2006} for PITA or \cite{ruprecht_contrib_ChenEtAl2014,ruprecht_contrib_DaiEtAl2013,ruprecht_contrib_GanderPetcu2008,ruprecht_contrib_RuprechtKrause2012} for Parareal, but so far none of these have been tested for the full Navier-Stokes equations.
\section{Conclusions}
The paper presents a numerical study of how the Reynolds number (or, inversely, the viscosity parameter) affects the convergence of the time-parallel Parareal method when used to solve the Navier-Stokes equations. 
From other works it is known that Parareal can develop a mild instability for problems with dominant imaginary eigenvalues, so it can be expected that as the viscosity is decreased, Parareal will eventually become unstable at some point.
A linear stability analysis is performed to motivate this assumption, which is then substantiated by numerical examples, solving a two-dimensional driven cavity problem for different Reynolds numbers and different spatial resolutions. 
It is confirmed that the convergence of Parareal deteriorates as the viscosity parameter becomes smaller and the flow becomes more and more dominated by convection.
This necessitates either the use of a very small time-step in the coarse method or many iterations of Parareal, but both these choices significantly reduce the achievable speedup.

\section*{Acknowledgments}
This work was supported by Swiss National Science Foundation (SNSF) grant 145271 under the lead agency agreement through the project ``ExaSolvers'' within the Priority Programme 1648 ``Software for Exascale Computing'' (SPPEXA) of the Deutsche Forschungsgemeinschaft (DFG).


%
%
%
\bibliographystyle{vmams}  

\input{ruprecht_contrib.bbl}

\end{document}

%% file: ruprecht_contrib.bbl
\ifx\undefined\bysame
\newcommand{\bysame}{\leavevmode\hbox to3em{\hrulefill}\,}
\fi